
\magnification=1200
\documentstyle{amsppt}
\hsize=152truemm
\vsize=215truemm
\hoffset=0truecm
\voffset=0,5truecm
\parindent=12pt
\parskip=0pt

\def\Int{{\Bbb Z}}

\def\bar{\overline}
\def\map{\longrightarrow}

\def\a{\alpha}
\def\al{\alpha}

\def\GL{\operatorname{GL}}
\def\PGL{\operatorname{PGL}}
\def\SL{\operatorname{SL}}

\def\Ep{\operatorname{Ep}}

\def\EU{\operatorname{EU}}
\def\FU{\operatorname{FU}}
\def\CU{\operatorname{CU}}

\def\Hom{\operatorname{Hom}}
\def\Max{\operatorname{Max}}
\def\Ann{\operatorname{Ann}}
\def\Rad{\operatorname{Rad}}

\def\Nat{{\Bbb N}}
\def\Co{{\Bbb C}}

\def\A{\operatorname{A}}

\def\F{\operatorname{F}}
\def\G{\operatorname{G}}
\def\K{\operatorname{K}}

\def\rk{\operatorname{rk}}

\def\unlhd{\trianglelefteq}
\def\lsup#1{{}^{#1}\!} 
\def\aeq{\Longleftrightarrow}
\def\seq{\Longrightarrow}
\def\map{\longrightarrow}

\def\ma{{\frak a}}




\magnification=1200




\def\Hom{\text{\rm Hom}}

\def\Ker{\text{\rm Ker}}

\def\rk{\operatorname{rk}}

\def\map{\longrightarrow}

\def\a{\alpha}

\def\bar{\overline}
\def\tilde{\widetilde}
\def\A{\operatorname{A}}

\def\F{\operatorname{F}}
\def\G{\operatorname{G}}
\def\K{\operatorname{K}}

\def\Ep{\operatorname{Ep}}

\def\SL{\operatorname{SL}}
\def\GL{\operatorname{GL}}

\def\GU{\operatorname{GU}}

\def\Max{\operatorname{Max}}

\def\sr{\operatorname{sr}}
\def\Int{{\Bbb Z}}

\vsize=215truemm
\hsize=152truemm

\topmatter
\title {The yoga of commutators}
\endtitle
\rightheadtext{the yoga of commutators}
\author{R.~Hazrat, A.~Stepanov, N.~Vavilov, Z.~Zhang}
\endauthor
\keywords
Unitary groups, Chevalley groups, elementary subgroups,
elementary generators, localisation, relative subgroups,
conjugation calculus, commutator calculus, Noetherian reduction,
Quillen---Suslin lemma, localisation-completion, commutator
formulae, commutator width, nilpotency of $\K_1$, nilpotent
filtration.
\endkeywords
\thanks
The first author 
acknowledges the support of EPSRC 
first grant scheme EP/D03695X/1.
The second author worked within the framework of the RFFI/Indian
Academy cooperation project 10-01-92651 ``Higher composition laws,
algebraic $K$-theory and algebraic groups'' (SPbGU--Tata Institute).
The third author worked within the framework of the RFFI/BRFFI
cooperation  project 10-01-90016 ``The structure of forms of
reductive groups, and behaviour of
small unipotent elements in representations of algebraic groups''
(SPbGU--Mathematics Institute of the Belorussian Academy).
Currently the work of the second and the third authors is supported
by the RFFI research projects 09-01-00762 (Siberian Federal
University), 09-01-00784 (POMI), 09-01-00878 (SPbGU),
RFFI/DFG cooperation project 09-01-91333 (POMI--LMU)
and RFFI/VAN coopera\-tion project 09-01-90304 (SPbGU--HCMU).
The fourth author acknowledges the support of NSFC grant 10971011
and the support from Beijing Institute of Technology.
\endthanks
\endtopmatter

\document

In the present paper we briefly describe three recent versions
of localisation me\-thods in the study of algebraic-like groups,
namely,
\smallskip
$\bullet$ {\bf Relative localisation} [24] -- [26]
\smallskip
$\bullet$ {\bf Universal localisation} [47],
\smallskip
$\bullet$ Enhanced {\bf localisation-completion}, [5], [17], [18],
[22], [7],
\smallskip\noindent
and state some recent results obtained therewith.
\par
The term {\tt yoga} {\tt of} {\tt commutators} refers to the methods
themselves, more precise\-ly, to a large body of calculations and
technical facts, conventionally known as the {\bf conjugation
calculus} and the {\bf commutator calculus}.
\par
As a matter of fact, the three typical recent applications of these
methods, we mention here, also pertain to commutators in
algebraic-like groups:
\smallskip
$\bullet$ Standard commutator formulae for congruence
subgroups/relative elementary subgroups, [67], [69], [24] -- [26].
\smallskip
$\bullet$ Universal bound for the width of arbitrary commutators
in terms of elementary generators [46], [49], [21], [47].
\smallskip
$\bullet$ Nilpotent structure of $\K_1$, for groups over a ring
of finite Bass---Serre dimension $d=\delta(R)<\infty$ it can be
interpreted as a multiple commutator formula of length $d+1$, see
[5], [17], [18], [22], [7].
\smallskip
This paper is based on our joint talks at the following Conferences:
\smallskip
$\bullet$ Topology, Geometry  and Dynamics: Rokhlin Memorial
(SPb, January 2010) [20],
\smallskip
$\bullet$ 2nd Group Theory Conference (Mashhad, Iran, March 2010),
\smallskip
$\bullet$ Polynomial Computer Algebra 2010 (SPb, April 2010),
\smallskip
$\bullet$ International Algebra Conference dedicated to the
70th Birthday of A.~V.~Ya\-kov\-lev (SPb, June 2010).
\smallskip\noindent
There, we described a major project whose goal is to review
existing localisation methods in the study of groups of
points of reductive algebraic groups, classical groups, and
related groups. Our main objective is to develop new more
powerful and efficient versions of conjugation calculus and
commutator calculus, with a view towards new applications.
\par
In this sense, it is a partial update of our survey [65], which
was based on our talk at the
\smallskip
$\bullet$ Applications of Computer Algebra 2008 (Linz, July 2008),
\smallskip
$\bullet$ Symbolic and Numeric Scientific
Computations 2008 (Linz, July 2008),
\smallskip
$\bullet$ Polynomial Computer Algebra 2009 (SPb, April 2009),
\smallskip\noindent
and had more computational flavour. At that time, we were rather
sceptical about the use of localisation techniques for actual
calculations in the groups of points of algebraic groups, with
realistic bounds.
\par
However, in the Fall and Winter 2009/2010 we developed new versions
of locali\-sa\-tion, which allowed us to prove some striking results.
For example, it turned out, that for algebraic groups length
estimates of various classes of elements, such as commutators,
in terms of elementary generators, do not depend on the
dimension of the ground ring, but on the type of the group alone.
\par
So far we have not succeeded in getting reasonable polynomial bounds,
even less so in converting our methods into working algorithms for
calculations in algebraic groups. Still, presently these goals seem
slightly less unfeasible, than at the moment we were writing [65].


\heading
\S~1. The groups
\endheading

Here, we consider algebraic-like or classical-like group
functors $G$. Further, let $G(R)$ be the group of points
of $G$ over a ring $R$. Observe, that groups of types
other than $\A_l$ only exist over commutative rings.
Typically, $G(R)$ is one of the following groups.
\smallskip
{\bf A.} \ General linear group $\GL(n,R)$ of degree $n$
over $R$.
\smallskip
Actually, many results are already new in this context.
Moreover, one can consider general linear groups over
arbitrary associative rings, and in this case our methods
work over quasi-finite rings. Recall, that a ring $R$
is called {\it almost commutative\/} (or, sometimes,
{\it module finite\/}), if it is finitely generated as a
module over its centre. {\it Quasi-finite} rings are
direct limits of inductive systems of almost commutative rings.
\par
However, in most of our current papers we work in one of
the following more general situations:
\smallskip
{\bf B.} \ Bak's unitary groups $\GU(2n,A,\Lambda)$, over a
form ring $(A,\Lambda)$.
\smallskip
The notation we use for these groups, their subgroups and elements
are mostly standard. As in [9], in the case of hyperbolic unitary
groups we number columns and rows of matrices as follows:
$1,\ldots,n,-n,\ldots,-1$.
Recall, that in this setting $A$ is a [not necessarily commutative]
ring with involution $\bar{\ }:A\map A$, $\Lambda$ is the form
parameter. To somewhat simplify matters, we usually assume that
$A$ is module finite over a commutative ring $R$. In general,
$\Lambda$ is not an $R$-module. Thus, $R$ has to be replaced
by its subring $R_0$, generated by $\xi\bar\xi$, for $\xi\in R$.
In the sequel we usually do not discuss similar technical details,
referring to [7], [9], [17], [18], [21], [23], [24], [28] for
precise statements and conclusive proofs.
\par
Actually, our favourite setting in this paper is the following
one, see [1] -- [3], [58], [64] and references there.
\smallskip
{\bf C.} \ Chevalley groups $G(\Phi,R)$ of type $\Phi$ over $R$.
\smallskip
Chevalley groups are indeed {\it algebraic\/}, and the ground
rings are {\it commutative\/} in this case, which usually makes
life easier. We illustrate most of our methods in this example.
\par
Together with the algebraic-like group $G(R)$ we consider the
following subgroups.
\smallskip
$\bullet$ First of all, the elementary group $E(R)$,
generated by elementary unipotents.
\smallskip
In the linear case, the elementary generators are elementary
[linear] transvections $t_{ij}(\xi)$, $1\le i\neq j\le n$, $\xi\in R$.
In the unitary case, the elementary generators are elementary
unitary transvections $T_{ij}(\xi)$, $1\le i\neq j\le -1$, $\xi\in A$.
In the even hyperbolic case they come in two modifications. They can
be short root type, $i\neq\pm j$, when the parameter $\xi$ can be
any element of $A$. On the other hand, for the long root type $i=-j$
and the parameter $\xi$ must belong to [something defined in terms of]
the form parameter $\Lambda$. Finally, for Chevalley groups,
the elementary generators are the elementary root unipotents
$x_{\a}(\xi)$ for a root $\a\in\Phi$ and a ring element $\xi\in R$.
\par
Further, let $I\unlhd R$ be an ideal of $R$. We also consider
the following relative subgroups.
\smallskip
$\bullet$ The elementary group $E(I)$ of level $I$, generated by
elementary unipotents of level $I$.
\smallskip
$\bullet$ The relative elementary group $E(R,I)=E(I)^{E(R)}$
of level $I$.
\smallskip
$\bullet$ The principal congruence subgroups $G(R,I)$ of level $I$,
the kernel of reduction homomorphism $\rho_I:G(R)\map G(R/I)$.
\smallskip
$\bullet$ The full congruence subgroups $C(R,I)$ of level $I$, the
inverse image of the centre of $G(R/I)$ with respect to $\rho_I$.
\smallskip
Recall the usual notation for these groups in the above contexts
A--C.
$$ \matrix
G(R)\qquad\qquad\hfill&\GL(n,R)\qquad\hfill&\GU(n,R,\Lambda)\qquad\hfill&G(\Phi,R)\hfill\\
\noalign{\vskip 5truept}
E(R)\qquad\qquad\hfill&E(n,R)\qquad\hfill&\EU(n,R,\Lambda)\qquad\hfill&E(\Phi,R)\hfill\\
\noalign{\vskip 5truept}
E(I)\qquad\qquad\hfill&E(n,I)\qquad\hfill&\FU(n,I,\Gamma)\qquad\hfill&E(\Phi,I)\hfill\\
\noalign{\vskip 5truept}
E(R,I)\qquad\qquad\hfill&E(n,R,I)\qquad\hfill&\EU(n,I,\Gamma)\qquad\hfill&E(\Phi,R,I)\hfill\\
\noalign{\vskip 5truept}
G(R,I)\qquad\qquad\hfill&\GL(n,R,I)\qquad\hfill&\GU(n,I,\Gamma)\qquad\hfill&G(\Phi,R,I)\hfill\\
\noalign{\vskip 5truept}
C(R,I)\qquad\qquad\hfill&C(n,R,I)\qquad\hfill&\CU(n,I,\Gamma)\qquad\hfill&C(\Phi,R,I)\hfill\\
\endmatrix $$
\par
There are two more general contexts, where Quillen---Suslin
localisation method has been successfully used by Victor Petrov,
Anastasia Stavrova, and Alexander Luzgarev,
[42] -- [44], [45], [36].
\smallskip
{\bf D.} \ Isotropic reductive groups $G(R)$,
\smallskip
{\bf E.} \ Odd unitary groups $U(V,q)$.
\smallskip\noindent
We are positive that one could obtain results similar to the ones
stated in the present paper also in these contexts, and we are
presently working towards it.


\heading
\S~2. Localisation
\endheading

In the present paper we only use commutative localisation. First,
let us fix some notation. Let $R$ be a commutative ring with
1, $S$ be a multiplica\-tive system in $R$ and $S^{-1}R$ be the
corresponding localisation. We mostly use localisation
with respect to the two following types of multiplicative
systems.
\par\smallskip
$\bullet$ {\it Principal localisation\/}: $S$ coincides with
$\langle s\rangle=\{1,s,s^2,\ldots\}$, for some non-nilpotent
$s\in R$, in this case we usually write $\langle s\rangle^{-1}R=R_s$.
\par\smallskip
$\bullet$ {\it Localisation at a maximal ideal}: $S=R\setminus M$,
for some maximal ideal $M\in\Max(R)$ in $R$, in this case
we usually write $(R\setminus M)^{-1}R=R_M$.
\par\smallskip
We denote by $F_S:R\map S^{-1}R$ the canonical ring homomorphism
called the {\it localisation homomorphism\/}. For the two special
cases above, we write $F_s:R\map R_s$ and $F_M:R\map R_M$,
respectively.
\par
When we write an element as a fraction, like $a/s$ or
$\displaystyle{\frac{a}{s}}$ we {\it always\/} think of it
as an element of some localisation $S^{-1}R$, where $s\in S$.
If $s$ were actually invertible in $R$, we would have written
$as^{-1}$ instead.
\par
Ideologically, all proofs using localisations are based on the
interplay of the three following observations:
\smallskip
$\bullet$ Functors of points $R\rightsquigarrow G(R)$ are
compatible with localisation,
$$ g\in G(R)\qquad\aeq\qquad F_M(g)\in G(R_M),
\quad\text{for all\ } M\in\Max(R). $$
\par
$\bullet$ Elementary subfunctors $R\rightsquigarrow E(R)$
are compatible with factorisation, for any $I\trianglelefteq R$
the reduction homomorphism $\rho_I:E(R)\map E(R/I)$ is surjective.
\smallskip
$\bullet$ On a [semi-]local ring $R$ the values of semi-simple
groups and their elementary subfunctors coincide, $G(R)=E(R)$.
\smallskip
The following property of the functors $G$ and $E$,
will be crucial for what follows: they {\it commute with direct
limits\/}. In other words, if $R=\varinjlim R_i$, where
$\{R_i\}_{i\in I}$ is an inductive system of rings, then
$$ X(\Phi,\varinjlim R_i)=\varinjlim X(\Phi,R_i). $$
\noindent
We use this property in the two following situations.
\par\smallskip
$\bullet$ {\it Noetherian reduction\/}: let $R_i$ be the inductive
system of all finitely generated subrings of $R$ with respect to
inclusion. Then
$$ X=\varinjlim X(\Phi,R_i). $$
\noindent
This allows to reduce most of the proofs to the case of
Noetherian rings.
\par\smallskip
$\bullet$ {\it Reduction to principal localisations\/}:
let $S$ be a multiplicative system in $R$ and let $R_s$, $s\in S$,
be the corresponding inductive system with respect to the
principal localisation homomorphisms: $F_{t}:R_s\map R_{st}$.
Then
$$ X(\Phi,S^{-1}R)=\varinjlim X(\Phi,R_s). $$
\noindent
This reduces localisation in any multiplicative system to
principal localisations.


\heading
\S~3. Injectivity of localisation homomorphism
\endheading

Most localisation proofs rely on the injectivity of localisation
homomorphism $F_S$. As observed in the previous section, we can
only consider {\it principal\/} localisation homomorphisms $F_s$.
Of course, $F_s$ is injective when $s$ is regular. Thus,
localisation proofs are particularly easy for integral domains.
A large part of what follows are various devices to fight with
the presence of zero-divisors.
\par
When $s$ is a zero-divisor, $F_s$ is not injective on the group
$G(\Phi,R)$ itself. But its restrictions to appropriate congruence
subgroups often are. Here are two important typical cases,
Noetherian rings and semi-simple rings.
\proclaim
{Lemma~1} Suppose\/ $R$ is Noetherian and\/ $s\in R$. Then
there exists a natural number $k$ such that the homomorphism\/
$F_s:G(\Phi,R,s^kR)\map G(\Phi,R_s)$ is injective.
\endproclaim
\demo
{Proof}
The homomorphism $F_s:G(\Phi,R,s^kR)\map G(\Phi,R_s)$ is
in\-jec\-ti\-ve when\-ever $F_s:s^kR\map R_s$ is injective.
Let $\ma_i=\Ann_R(s^i)$ be the annihilator of $s^i$ in $R$.
Since $R$ is Noetherian, there exists $k$ such that
$\ma_k=\ma_{k+1}=\ldots$. If $s^ka$ vanishes in $R_s$, then
$s^is^ka=0$ for some $i$. But since $\ma_{k+i}=\ma_k$, already
$s^ka=0$ and thus $s^kR$ injects in $R_s$.
\enddemo
\proclaim
{Lemma~2}
If\/ $\Rad(R)=0$, then\/ $F_s:G(\Phi,R,sR)\map G(\Phi,R_s)$
is injective for all\/ $s\in R$, $s\neq 0$.
\endproclaim
\demo
{Proof}
It suffices to prove that $F_s:sR\map R_s$ is injective.
Suppose that $s\xi\in sR$ goes to $0$ in $R_s$. Then there
exists an $m\in\Nat$ such that $s^ms\xi=0$. It follows that
$(s\xi)^{m+1}=0$ and since $R$ is semi-simple, $s\xi=0$.
\enddemo
In [22] we used reduction to Noetherian rings, whereas in
[49] reduction to semi-simple rings was used.
\par
Another important trick to override the presence of zero-divisors
consists in throwing in polynomial variables. Namely, instead
of the ring $R$ itself we consider the polynomial ring $R[t]$
in the variable $t$. In that ring $t$ is not a zero-divisor,
so that the localisation homomorphism $F_t$ is injective.
We can use that, and then specialise $t$ to any $s\in R$.
\par
Actually, throwing in polynomial variables has mor� than one
use. The �lemen\-ta\-ry subfunctors $R\rightsquigarrow E(R)$ are
not compatible with localisation,
$$ g\in E(R)\qquad\seq\qquad F_M(g)\in E(R_M),
\quad\text{for all\ } M\in\Max(R), $$
\noindent
but the converse implication does not hold, for otherwise
$E(R)$ would coincide with [the semi-simple part of] $G(R)$
for all commutative rings.
\par
The following remarkable observation was due to Daniel Quillen
at the level of $\K_0$, and was first applied by Andrei Suslin
at the level of $\K_1$, in the context of solving Serre's
conjecture, and its higher analogues [51]. See [30] for a
description of Quillen---Suslin's idea in its historical
development. We refer to the following result as Quillen---Suslin's
lemma.
\proclaim
{Theorem 1} Let $g\in G(R[t],tR[t])$. Then,
$$ g\in E(R[t]) \qquad\aeq\qquad F_M(g)\in E(R_M[t]),
\quad\text{for all\ } M\in\Max(R).  $$
\endproclaim


\heading
\S~4. How localisation works
\endheading

As was already mentioned, {\bf localisation and patching}
was first used to study the structure of linear groups by
Andrei Suslin, back in 1976, see [51]. Among other
important early contributors one could mention Suslin's
[then] students Vyacheslav Kopeiko [29], [52] and Marat
Tulenbaev [55], as well as Leonid Vaserstein [56], [58] -- [60],
Eiichi Abe [1] -- [3], Li Fuan [31], [32], Giovanni Taddei
[53], [54], You Hong [62].
\par
Let us illustrate how localisation works in the classical
example, normality of the elementary subfunctor. Namely,
we wish to prove that $E(R)\trianglelefteq G(R)$, for any
commutative ring $R$.
\par
To be more specific, below we assume that $G(R)=G(\Phi,R)$
is the simply-connected Chevalley group of type $\Phi$. In
this case, $E(R)=E(\Phi,R)$ is generated by the elementary
root unipotents $x_{\a}(\xi)$, for $\a\in\Phi$ and $\xi\in R$.
From here on, we {\it always\/} assume that {\smc $\Phi$ is
reduced and irreducible of rank $\ge 2$}. We do not reproduce
this standing assumption in the statements of subsidiary
results.
\par
Thus, we wish to prove that for all $g\in G(\Phi,R)$, all
$\al\in\Phi$ and all $\xi\in R$ one has
$$ x=gx_{\al}(\xi)g^{-1}\in E(\Phi,R). $$
\noindent
All localisation proofs are based on {\it partitions of $1$\/}.
In other words, we pick up $\zeta_1,\ldots,\zeta_m\in R$ such that
$\zeta_1+\ldots+\zeta_m=1$ and each of $gx_{\al}(\zeta_i\xi)g^{-1}$
already lies in $E(\Phi,R)$. The difference between various
localisation methods is in how one chooses such a partition.
\par
For our taste, the most elementary way is the following version
of {\bf localisation and patching} method. Instead of throwing
in independent variables, as Quillen and Suslin did originally,
and many others after them, we follow Anthony Bak [5], and
recourse to Noetherian reduction.
As we observed, the functors $G=G(\Phi,\underline{\ \ }\ )$ and
$E=E(\Phi,\underline{\ \ }\ )$ commute with direct limits.
Since $R$ is the direct limit of its finitely generated subrings,
we can from the very start assume that $R$ is Noetherian.
This will allow us to invoke Lemma~1.
\par
Since we work with simply connected groups, for a local ring $R$
the elementary subgroup $E(\Phi,R)$ coincides with the Chevalley
group $G(\Phi,R)$. Thus, for any maximal ideal $M\in\Max(R)$
one has $F_M(g)\in E(\Phi,R_M)$. Now, we again invoke the fact
that the functors $G$ and $E$ commute with direct limits.
Since $R_M$ is the direct limit of $R_t$, $t\in R\setminus M$,
there exists an $s\in R\setminus M$ such that
$F_s(g)\in E(\Phi,R_s)$.
\par
We will search for $\zeta_i$'s in the desired partition of 1, as
multiples of high powers $s^l$ of the above elements $s$, for
various maximal ideals $M\in\Max(R)$ and sufficiently large
exponents $l$. Set $y=gx_{\al}(s^l\eta\xi)g^{-1}$, for some
$\eta\in R$. Since the ring $R$ is Noetherian, we can apply
Lemma~1, and conclude that for a large power of $s$, say for
$s^n$, the restriction of $F_s$ to the principal congruence
subgroup $G(\Phi,R,s^nR)$ is injective.
\par
First, we argue locally, this part of the proof is called
[first] {\bf localisation\/}. Since $F_s(g)\in E(\Phi,R_s)$, it
can be written as a product of elementary root unipotents
$x_{\a}\big(F_s({\theta}/{s^k})\big)$, $\theta\in R$, $k\ge 0$.
From the Chevalley commutator formula it follows that conjugation
by such an element is continuous in $s$-adic topologuy, this is
exactly the {\bf conjugation calculus} we discuss in the next
section. In particular, there exists a high power of $s$, say,
$s^l$, $l\gg n$, such that
$$ F_s(y)=F_s(g)F_t(x_{\a}(s^l\eta\xi))F_s(g)^{-1} $$
\noindent
can be expressed as a product
$$ F_s(y)=\prod_{j=1}^m x_{\beta_j}(F_s(s^n\theta_j))\in E(\Phi,R_s), $$
\noindent
for some $\theta_j\in R$.
\par
Take the product
$$ z=\prod_{i=1}^m x_{\beta_i}(s^nc_i)\in E(\Phi,R), $$
\noindent
By the very definition $F_s(z)=F_s(y)$. On the other hand, since
$G(\Phi,R,s^nR)$ is normal in $\G$, one has $y,z\in G(\Phi,R,s^nR)$.
Injectivity of $F_s$ implies that $y=z\in E(\Phi,R)$.
\par
The final part of the proof is called {\bf patching}. Since
$s^l\notin M$ and the same works for all maximal ideals, we can
choose a finite set of such powers $s_i^{l_i}$, $1\le i\le m$,
which generate $R$ as an ideal,
$$ \zeta_1+\ldots+\zeta_m=s_1^{l_1}\eta_1+\ldots+s_m^{l_m}\eta_m=1 $$
\noindent
is the desired partition.
\par
Of course, there are some further technical details. For example,
when one works with the adjoint group, such as $\PGL(n,R)$ or a
diagonal extension of a Chevalley group, such as $\GL(n,R)$,
there is an extra toral factor to take care of.


\heading
\S~5. Conjugation calculus
\endheading

The first main objective of the conjugation calculus is to establish
that conjuga\-tion by a fixed matrix $g\in G(\Phi,R_s)$ is continuous
in $s$-adic topology. In the proof one uses a base of $s$-adic
neighborhoods of $e$ and establishes that for any such neighborhood
$V$ there exists another neighborhood $U$ such that ${}^gU\subseteq V$.
\par
To be more specific, let us state some typical results of conjugation
calculus for Chevalley groups. Usually, as the base of neighborhoods
of $e$, one takes
\smallskip
$\bullet$ elementary subgroups $E(\Phi,s^mR)$, or
\smallskip
$\bullet$ relative elementary subgroups $E(\Phi,R,s^mR)$.
\smallskip
For advanced applications, one usually needs more than just continuity
of conju\-ga\-tion by $g$. One has to estimate the {\it module of
continuity\/}, depending on the size of denominators in expression of
$g$ as a product of $s$-elementary factors. In generation problems, one
often has to estimate also the {\it length\/} of arising elementary
expressions.
\par
To state typical results in this direction, we have to introduce some
further notation. Namely, let $L$ be a nonnegative integer and let
$E^L(\Phi,I)$ denote the {\it subset\/} of $E(\Phi,I)$ consisting of
all products of $L$ or fewer elementary root unipotents $x_{\a}(\xi)$,
where $\a\in\Phi$ and $\xi\in I$. Thus, $E^1(\Phi,I)$ is the set of
all $x_{\a}(\xi)$, $\a\in\Phi$, $\xi\in I$.
\par
Conjugation calculus and commutator calculus are rare examples of
induction results, where the base of induction is {\it terribly\/}
much harder, than the induction step. Without length estimates the
following results have been established by Giovanni Taddei [53], [54],
and then, in a stronger and more straightforward form, by the first
and the third author [22]. The precise form with explicit length
estimates is taken from the paper by the second and the third author
[49].
\proclaim
{Lemma~3} If\/ $p,q$ and\/ $h$ are given, there exist\/ $o,r$
such that
$$ \lsup xy\in E^{24}(\Phi,s^pt^qR),\qquad
\text{for all}\quad
x\in{E^1\Bigl(\Phi,\frac{1}{s^h}R\Bigr)},\quad
y\in{E^1\Bigl(\Phi,{s^{o}t^r}R\Bigr)}. $$
\endproclaim
For the case, where $x$ and $y$ are not opposite, the claim
immediately follows from the Chevalley commutator formula.
Indeed, let $i_{\Phi}$ be the largest integer which may appear as
$i$ in a root $i\a+j\beta\in\Phi$ for all $\a,\beta\in\Phi$.
Obviously $i_{\Phi}=1,2$ or 3, depending on whether $\Phi$ is
simply laced, doubly laced or triply laced.
\par
Now, let $\a\neq-\beta$ and set $o\geq i_{\Phi}h+p+1$,
$r\ge q$. By the Chevalley commutator formula, one has
$$ x_{\a}\Bigl(\frac{a}{s^h}\Bigr)
x_{\beta}\Bigl(s^ot^rb\Bigr)
x_{\a}\Bigl(\hskip-2pt-\frac{a}{s^h}\Bigr)=
x_{\beta}\Bigl(s^ot^rb\Bigr)\prod_{i\a+j\beta\in\Phi}x_{i\a+j\beta}
\biggl(N_{\a\beta ij}{\Bigl(\frac{a}{s^h}\Bigr)}^i{\Bigl(s^ot^rb\Bigr)}^j\biggr)
$$
\noindent
and a quick inspection shows that the right hand side of the
above equality is in $E^L(\Phi,s^pt^qR)$, where $L=2,3$ or $5$,
depending on whether $\Phi$ is simply laced, doubly laced or
triply laced.
\par
For the case of opposite roots, one first has to use the Chevalley
commutator formula to express $x_{-\a}\Bigl(s^ot^rb\Bigr)$ as
a product of elementary factors, corresponding to the roots not
opposite to $\a$. For example, when $-\a=\gamma+\delta$ is the sum
of two roots of the same length, one has
$$ x_{-\a}(s^ot^rb)=
\bigl[x_{\gamma}(s^{o/2}t^{r/2}),x_{\delta}(s^{o/2}t^{r/2}b)\bigr], $$
\noindent
and we have reduced the problem to the preceding case. For other
cases the proof is similar, but slightly fancier, due to the longer
products in the Chevalley commutator formula, see [22], [49] for
details.
\par
Now, the following general result immediately follows by induction.
\proclaim
{Lemma~4} If\/ $p,q$ and\/ $h$ are given, there exist\/ $o,r$
such that
$$ \lsup xy\in E^{24^LK}(\Phi,s^pt^qR),\qquad
\text{for all}\quad
x\in{E^L\Bigl(\Phi,\frac{1}{s^h}R\Bigr)},\quad
y\in{E^K\Bigl(\Phi,{s^{o}t^r}R\Bigr)}. $$
\endproclaim
Actually, for systems without factors of type $\G_2$ one can
even conclude that $\lsup xy\in E^{13^LK}(\Phi,s^pt^qR)$.
For simply laced systems and for $\F_4$, one can conclude
that $\lsup xy\in E^{8^LK}(\Phi,s^pt^qR)$.


\heading
\S~6. Commutator calculus
\endheading

More sophisticated applications, such as calculation of mutual
commutator sub\-groups, nilpotent filtration, description of various
classes of intermediate subgroups, etc., require {\bf second
localisation}. In other words, we have to be able to
simulta\-ne\-ously
fight with powers of {\it two\/} elements in the denominator.
\par
For $\GL(n,R)$ second localisation was used by Anthony Bak in [5],
and then generalised to unitary groups in the Thesis of the first
named author [17], [18]. For Chevalley groups, it was first 
implemented by the present authors in [22], and then in a more 
precise form in [49].
\proclaim
{Lemma~5} Given\/ $s,t\in R$ and\/ $p,q,k,m\in\Nat$ there
exist\/ $l,n\in\Nat$ and $L=L(\Phi)$ such that
$$ [x,y]\in E^{L}\bigl(\Phi,s^pt^qR\bigr),\qquad
\text{for all}\quad
x\in{E^1\Bigl(\Phi,\frac{t^{l}}{s^k}R\Bigr)},\quad
y\in{E^1\Bigl(\Phi,\frac{s^{n}}{t^m}R\Bigr)}. $$
\endproclaim
Let $\alpha,\beta\in\Phi$ and $a,b\in R$. We have
to prove that the commutator
$$ \biggl[x_\alpha\Bigl(\frac{t^l}{s^k}a\Bigr),
x_\beta\Bigl(\frac{s^n}{t^m}b\Bigr)\bigg]\in E^{L}\bigl(\Phi,s^pt^qR\bigr), $$
\noindent
for some specific $L$. For the case, where $\alpha\neq-\beta$ the
proof is easy. Writing the Chevalley commutator formula
$$ \biggl[x_\a\Bigl(\frac{t^l}{s^k}a\Bigr),
x_{\beta}\Bigl(\frac{s^{n}}{t^{m}}b\Bigr)\biggr]=
\prod_{i,j>0}x_{i\a+j\beta}\biggl(\Bigl(\frac{t^l}{s^k}a\Bigr)^i
\Bigl(\frac{s^{n}}{t^{m}}b\Bigr)^{\hskip-2pt j\ }\biggr), $$
\noindent
we see that one can take $l$ and $n$ large enough to kill
the denominators on the right hand side, and still leave
large enough powers of $s$ and $t$ in the numerators.
In fact, $l\ge i_{\Phi}m+q+1$ and $n\ge i_{\Phi}k+p+1$ would go.
Furthermore, the number of factors on the right hand side
of the Chevalley commutator formula is not more than 4. Thus,
the product on the right hand side is in $E^L(\Phi,s^pt^qR)$,
where $L=1,2$ or $4$, depending on whether $\Phi$ is simply
laced, doubly laced or triply laced.
\par
The proof for opposite roots is {\it much\/} fancier, and relies
on longer commutator identities.
\par
Again, the general case easily follows by induction, via purely
group theoretic arguments, see [49], \S~9.
\proclaim
{Lemma~6} Let\/ $s,t\in R$ and\/ $p,q,k,m\in\Nat$ Then there
exist\/ $l,n\in\Nat$ such that
$$ [x,y]\in E^{(L+1)^K-1}(\Phi,s^pt^qR)\qquad
\text{for all }
x\in E^1\Big(\Phi,\frac{s^{l}}{t^k}R\Big),\quad
y\in E^K\Bigl(\Phi,\frac{t^{n}}{s^m}R\Bigr). $$
\endproclaim
By looking inside the proofs in [22] and [49] one gets the
following silly length estimates
\smallskip
$\bullet$ $L\le 585$, for simply laced systems,
\smallskip
$\bullet$ $L\le 61882$, for doubly laced systems,
\smallskip
$\bullet$ $L\le 797647204$, for triply laced systems.
\smallskip\noindent
In obtaining these stupid bounds we do not look inside the
commutators and do not count the actual factors appearing in
the Chevalley commutator formula. However, should we do that,
the resulting bound for $\G_2$ still would be well in the millions.
\proclaim
{Problem~1} Calculate realistic length bounds for\/ $L$ in
these lemmas.
\endproclaim
As far as we can see, calculating in the 7-dimensional or the
8-dimensional representa\-tion of the group of type $\G_2$ one
gets bounds for $L$ within few dozens, rather than millions.
\par
Let us state a typical target result of the commutator calculus.
\proclaim
{Theorem 2} Fix an element $s\in R$, $s\neq 0$. Then for any $p$
and $k$ there exists an $r$ such that
$$ \Big[E\big(\frac{1}{s^k}R\big),F_s(G(\Phi,R,s^rR))\Big]\le
E(\Phi,F_s(s^pR))\le G(\Phi,R_s). $$
\endproclaim
Despite its rather technical appearance, it is a very general
and powerful result. In fact, in the {\it trivial\/} special
case, where $s=1$, this Theorem boils down to the normality of
the elementary subgroup!
\par
The general case of the theorem was substantially used in
description of over\-groups of classical and exceptional groups by
the second named author and Victor Petrov [66], and
by Alexander Luzgarev [35]. We do not mention any further results in
this direction, referring to our surveys [68] and [50].


\heading
\S~7. Relative commutator calculus
\endheading

For the group $G(R)$ itself, conjugation calculus works marvelously,
as one takes $E(s^mR)$ or $E(R,s^mR)$, as the base of $s$-adic
neighbourhoods. But can one {\it relativise\/} all occurring
calculations? In other words, what happens when we replace the
ring $R$ by an ideal $I\trianglelefteq R$? Again, one has to
establish that for any neighbourhood $V$ of $e$ in $G(R,I)$
there exists another neighborhood $U$ such that ${}^gU\subseteq V$.
\par
As a first attempt, without much thinking, one tries to replace
$R$ by $I$ everywhere in the above calculations. For example, it
seems that one should consider the following bases of $s$-adic
neighborhoods of $e$ in $G(R,I)$:
\smallskip
$\bullet$ elementary subgroups $E(s^mI)$,
\smallskip
$\bullet$ relative elementary subgroups $E(R,s^mI)$.
\smallskip
\par
However, both choices are not fully satisfactory in that they lead to
extremelly onerous calculations. The reason is that the first of these
choices is too small as the neighbourhood on the right hand side,
while the second of these choices is too large as the neighbourhood
on the left hand side.
\par
Solving problems posed by two of the present authors in [67], the
first and the last authors proposed in [26] a first fully functional
version of localisation at the relative level. Their idea was to take
the following {\it partially\/} relativised base of $s$-adic
neighbourhoods.
$$ E(s^mR,s^mI)={E(s^mI)}^{E(s^mR)}. $$
\par
To convey the flavour of the ensuing results, let us state some typical
lemmas from our forthcoming relative Chevalley paper [25]. Similar
results for $\GL(n,R)$ and $\GU(2n,R,\Lambda)$ are established in
[26] and [24], respectively. Again, the base of induction is the
hardest part of the whole argument.
\proclaim
{Lemma~7}
If $p$, $q$ and $k$ are given, there exist $h$ and $m$ such that
$$ {}^{\textstyle{E^1\left(\Phi,\frac{1}{s^k}R\right)}}
E(\Phi,s^ht^m I)\subseteq
E(\Phi,s^pt^qR,s^pt^q I). $$
\noindent
Such $h$ and $m$ depend on $\Phi$, $k$, $p$ and $q$ alone, but not on
the ideal $I$.
\endproclaim
Observe, that the proofs of this and similar results work in terms of
roots alone, and thus one obtains {\it uniform\/} estimates for the
powers of $s$ and $t$, which do not depend on the ideal $I$. In other
words, conjugation by $g\in G(\Phi,R_s)$ is {\it equi-continuous\/}
in all congruence subgroups $G(\Phi,R,I)$, with respect to
$$ E(\Phi,s^kR,s^k I)={E(\Phi,s^k I)}^{E(\Phi,s^kR)}, $$
\noindent
as the corresponding bases of $s$-adic neighbourhood.
\par
This is extremely important for applications we have in mind. For
example, in the next result we use this to obtain a uniform bound
for {\it two\/} ideals $A,B\trianglelefteq R$.
\proclaim
{Lemma~8}
If $p,k$ are given, then there is an $q$ such that
$$ {}^{E^1\left(\Phi,\frac{R}{s^k}\right)} \big
[E(\Phi,s^qR,s^q A), E(\Phi,s^qR,s^q B)\big] \subseteq
\big [E(\Phi,s^pR,s^p A), E(\Phi,s^pR,s^p B)\big]. $$
\endproclaim
Similarly, the induction base of the relative commutator calculus
looks as follows. Again, in view of applications, we state it for
{\it two\/} ideals $A,B\trianglelefteq R$.
\proclaim
{Lemma~9} If $p,q,k,m$ are given, then there exist $l$ and $n$
such that
$$ \Big [E^1\Big(\Phi,\frac{t^l}{s^k} A\Big),
E^1\Big(\Phi,\frac{s^n}{t^m}B\Big)\Big]
\subseteq \big[E(\Phi,s^pt^qR,s^pt^q A),E(\Phi,s^pt^qR,s^pt^q B)]. $$
\noindent
These $l$ and $n$ depend on $\Phi,p,q,k,m$ alone, and do not depend
on the choice of ideals $A$ and $B$.
\endproclaim
This is a rather difficult technical result. Also, in the relative
case induction step itself is non-trivial. For example, Lemma 9 itself
does not suffice even to start the induction. Instead, we have to
establish something as follows:
$$ \bigg[E^1(\Phi,s^q A),{}^{E^1\left(\Phi,\frac{R}{s^k}\right)}
E^1\Big(\Phi,\frac{B}{s^k}\Big)\bigg]
\subseteq \big[E(\Phi,s^pR,s^p A),E(\Phi,s^pR,s^p B)\big]. $$
\par
We do not reproduce the precise target results of the relative
commutator calcu\-lus, which are far too technical for a casual
overview. The interested reader can find such precise statements
in our papers [21], [24] -- [26].


\heading
\S~8. Relative commutator formulae
\endheading

One of the central and most important results in the theory of linear
groups over rings are the {\it absolute\/} standard commutator
formulae
$$ [G(R),E(R,I)]=E(R,I)=[E(R),C(R,I)]. $$
\noindent
The first one of them amounts to saying that $E(R,I)$ is normal in
$G(R)$, while the second one is a key tool of level reduction.
For $\GL(n,R)$ at the stable level these formulae were established
by Hyman Bass. Later, Leonid Vaserstein improved the estimate for
the stable rank $\sr(R)$ of the ring $R$ by 1.
\par
Soon thereafter, Andrei Suslin [51], Leonid Vaserstein [56],
Zenon Borewicz and the first author, observed that for {\it almost
commutative\/} rings these formulae hold in $\GL(n,R)$ for any $n\ge 3$.
Recall, that a ring $R$ is called almost commutative if it is finitely
generated as a module over its centre.
\par
In the sequel, Igor Golubchik, Alexander Mikhalev sen., Sergei Khlebutin
and Anthony Bak generalised these formulae to broad classes of
non-commutative rings. Vyacheslav Kopeiko and Andrei Suslin [29],
[52], Giovanni Taddei, Leonid Vaser\-stein, You Hong, Anthony Bak, and
the present authors generalised these results to other groups.
In [6], [8] -- [10], [16], [19], [22], [23], [30], [48], [64], [70]
one can find different proofs of these results, {\it many\/} further
references, and a detailed discussion of their role in the structure
theory.
\par
However, much less was known about the {\it relative\/} versions
of the above formulae. Namely, let $A,B\trianglelefteq R$ be two ideals
of the ring $R$. What can be said about the mutual commutators of
congruence subgroups and relative elementary subgroups of levels $A$
and $B$?
\par
Before our works this problem was only addressed at the stable level,
by Alec Mason and Wilson Stothers [37] -- [40].
\par
Let us state some typical results, we can prove by relative
versions of localisation methods, described in the previous section,
see [24] -- [26], [67], [69].
\proclaim
{Theorem 3A} Let\/ $R$ be a quasi-finite ring,\/ $n\ge 3$. Then for any
two ideals\/ $A,B\trianglelefteq R$ one has
$$ [E(n,R,A),\GL(n,R,B)]=[E(n,R,A),E(n,R,B)]. $$
\endproclaim
\proclaim
{Theorem~3B}
Let\/ $n\ge 3$,\/ $R$ be a commutative ring,\/ $(A,\Lambda)$ be a form
ring such that\/ $A$ is a quasi-finite\/ $R$-algebra. Further, let\/
$(I,\Gamma)$ and\/ $(J,\Delta)$ be two form ideals of a form ring\/
$(A,\Lambda)$. Then
$$ \big[\EU(2n,I,\Gamma),\GU(2n,J,\Delta)\big]=
\big[\EU(2n,I,\Gamma),\EU(2n,J,\Delta)\big]. $$
\endproclaim
\proclaim
{Theorem~3C}
Let\/ $\Phi$ be a reduced irreducible root system,\/ $rk(\Phi)\ge 2$.
Further, let\/ $R$ be a commutative ring, and\/ $A,B\trianglelefteq R$ be
two ideals of\/ $R$. Then
$$ [E(\Phi,R,A),G(\Phi,R,B)]=[E(\Phi,R,A),E(\Phi,R,B)]. $$
\endproclaim
Observe, that in general one cannot expect the equality
$$ [E(R,A),E(R,B)]=E(R,AB). $$
\noindent
However, the true reason, why this equality holds in the absolute
case, is not the fact that one of the ideals $A$ or $B$ coincides
with $R$, but just the fact that $A$ and $B$ are comaximal. Namely,
by combining the above relative commutator formulae with commutator
identities, like the celebrated Hall---Witt identity, one gets the
following results.
\proclaim
{Theorem~4A} Let\/ $R$ be a quasi-finite ring,\/ $n\ge 3$. Then for any
two comaximal ideals\/ $A,B\trianglelefteq R$, $A+B=R$, one has
$$ [E(n,R,A),E(n,R,B)]=E(n,R,AB+BA). $$
\endproclaim
We do not recall notation pertaining to form ideals [9], [14] --
[18], [23], [24].
\proclaim
{Theorem~4B} Let\/ $n\ge 3$, and\/ $(A,\Lambda)$ be an arbitrary form
ring for which absolute standard commutator formulae are satisfied.
Then for any two comaximal form ideals\/ $(I,\Gamma)$ and $(J,\Delta)$
of the form ring $(A,\Lambda)$, $I+J=A$, one has the following equality
$$ [\EU(2n,I,\Gamma),\EU(2n,J,\Delta)]=
\EU(2n,IJ+JI,{}^J\Gamma+{}^I\Delta+\Gamma_{\min}(IJ+JI)). $$
\endproclaim
As for the next result, we also have a more general version, in terms
of admissible pairs [2], [3], [16], [64]. We do not reproduce it here,
not to overburden the reader with technical details.
\proclaim
{Theorem~4C}
Let $\Phi$ be a reduced irreducible root system, $\rk(\Phi)\ge 2$.
Further, let $R$ be a commutative ring, and $A,B\trianglelefteq R$ be
two comaximal ideals of $R$, $A+B=R$, one has the following equality
$$ [E(\Phi,R,A),E(\Phi,R,B)]=E(\Phi,R,AB). $$
\endproclaim


\heading
\S~9. Anti-Ore
\endheading

Over fields, the groups of points of algebraic groups essentially
consist of com\-mu\-ta\-tors:
In fact, the celebrated {\it Ore conjecture\/} -- now a theorem,
[13], [34] -- asserts that every element of a [non-abelian]
finite simple group is a single commutator. It is usually only
easier to establish similar results for infinite fields.
\par
The results we formulate in this and the next sections go
in the opposite direction. Morally, they say that
{\smc groups of points of algebraic groups over rings have very
few commutators}. In a strict technical sense, they have not more
commutators, than elementary generators!
\par
Similar bounded width results have a long history, which we cannot
even sketch here.
In general, $G(R)$ does not have bounded width with respect to the
elementary generators.
\smallskip
$\bullet$ First, it is not even spanned by them! By definition,
elementary generators generate the {\it elementary\/}
subgroup $E(R)$, which is usually strictly smaller, than $G(R)$.
\smallskip
$\bullet$ Even when $G(R)=E(R)$, it does not have to have bounded
width with respect to the elementary generators. Wilberd van der
Kallen observed that already $\SL(3,\Co[x])$ has unbounded
width [27].
\smallskip
For some time, it was an open question, whether $E(R)$ has bounded
length with respect to commutators. It was settled in the negative
by Keith Dennis and Leonid Vaserstein [11], [12].
\par
However, the situation with the commutators turned out to be
exactly the opposite to what was expected in the 1980-ies.
The following amazing result is established in the paper
by Alexander Sivatsky and the second named author [46].
\proclaim
{Theorem~5A} Let\/ $G=\GL(n,R)$, $n\ge 3$, where\/ $R$ be a
Noetherian ring such that\/ $\dim\Max(R)=d<\infty$. Then there
exists a natural number\/ $N$ depending only on\/ $n$ and\/ $d$
such that each commutator\/ $[x,y]$ of elements\/
$x\in\GL(n,R)$ and\/ $y\in E(n,R)$ is a product of at most\/ $N$
elementary transvections.
\endproclaim
The original proof of that result in [46] depended {\it both\/}
on localisation and a very precise form of decomposition of
unipotents, as proven in [48]. It was not at all clear, that it
could be generalised to other groups, even the classical ones.
\par
However, the second named author soon came up with an idea to
replace the use of decomposition of unipotents by the second
localisation. Using the machinery developed by the first named
and the third named authors in [22] --- and, in fact, enhancing
it --- the second and the third named authors succeded in
generalising the above result to Chevalley groups [49].
\proclaim
{Theorem~5C} Let\/ $G=G(\Phi,R)$ be a Chevalley group of rank\/
$l\ge 2$ and let\/ $R$ be a ring such that\/ $\dim\Max(R)=d<\infty$.
Then there exists a natural number\/ $N$ depending only on\/ $\Phi$
and\/ $d$ such that each commutator\/ $[x,y]$ of
elements\/ $x\in G(\Phi,R)$ and\/ $y\in E(\Phi,R)$ is a product
of at most\/ $N$ elementary root unipotents.
\endproclaim
We do not describe the strategy of the localisation proof of this
result, since it is expounded in full detail in [49]. Moreover,
this result is largely superceded by the truly miraculous result
stated in the next section.
\par
Similar result also holds for unitary groups [21], it relies on
the full force of localisation methods developed in [17], [18], [25].
\proclaim
{Theorem~5B}
Let\/ $G=\GU(2n,R,\Lambda)$ be the unitary group of degree\/ $2n$
over a commutative form ring $(R,\Lambda)$. Assume that\/ $n\ge 2$
and\/ $\dim\Max(R)=d<\infty$. Then there exists a natural number\/ $N$
depending only on\/ $n$ and\/ $d$ such that each commutator\/ $[x,y]$
of elements\/ $x\in\GU(2n,R,\Lambda)$ and\/ $y\in\EU(2n,R,\Lambda)$
is a product of at most\/ $N$ elementary unitary transvections.
\endproclaim


\heading
\S~10. Universal localisation
\endheading

Now, something truly amazing will happen. The results we stated
in the previous section relied on the fact that Jacobson
dimension $d=\dim(\Max(R))$ of the ground ring $R$ is finite.
The length estimates of commutators, in elementary generators,
were stated in terms of degree $n$ or type of root system $\Phi$,
{\it and\/} dimension $d$.
\par
Recently, the second named author observed, that for {\it algebraic\/}
groups no finite\-ness condition on the ground rings is necessary
here. In other words, the length estimates do not depend on $d$.
In particular, {\smc there exist universal length bounds for
the length of commutators}, in the group of given type, over an
{\it arbitrary\/} commutative ring. One does not even have to
assume these rings to be Noetherian!
\par
Let us state one of the main results of [47].
\proclaim
{Theorem~6C} Let $\Phi$ be a reduced irreducible roots system of
rank $l\ge 2$, and let $G=G(\Phi,\underline{\ \ }\ )$ be the simply
connected Chevalley---Demazure group scheme of type $\Phi$. Then
there exists an integer $l$ depending only on $\Phi$, that
satisfies the following property. For any commutative ring $R$,
any $x\in G(\Phi,R)$ and any $y\in E(\Phi,R)$ the commutator
$[x,y]$ can be written as a product of at most $l$ elementary
root unipotents in $G(\Phi,R)$.
\endproclaim
The idea behind this result can be described as follows. There
exists a {\bf universal commutator} which is {\it generic\/}
in the sense that it specialises to any other commuta\-tor of this
shape. Thus, any elementary expression of this universal
commutator provides an upper bound on the length of any such
commutator. One expects such a universal commutator to live inside
the group of points over the {\bf universal coefficient ring}
of the group.
\par
It is easy to guess, what is the universal coefficient ring for
the algebraic group itself. Recall, that
$$ G(\Phi,R)=\Hom(\Int[G],R), $$
\noindent
where $\Int[G]$ is the affine ring of $G=G(\Phi,\underline{\ \ }\ )$.
In other words, a point $h\in G(\Phi,R)$ of the group
$G(\Phi,\underline{\ \ }\ )$ over the ring $R$ can be identified
with a homomorphism $h:\Int[G]\map R$.
\par
Clearly, any such homomorphism $h$ can be factored through the
identity map $\Int[G]\map\Int[G]$. In other words, this means
that the point $g\in G(\Phi,\Int[G])$ represen\-ted by the
identity map $\Int[G]\map\Int[G]$, is the {\bf generic element}
of $G$ in the sense that it specialises to any point
$h\in G(\Phi,R)$ over any commutative ring $R$. Namely, by
functoriality,
$$ G(h):G(\Phi,\Int[G])\map G(\Phi,R),\quad g\mapsto h. $$
\par
Similarly, one can define {\bf generic elements} under
localisation. Namely, for $s\in\Int[G]$ we denote by
$g_s\in G(\Phi,\Int[G]_s)$ the element, represented by the
localisation homomorphism $F_s:\Int[G]\map \Int[G]_s$.
\par
As a special case of this construction, we get generic
{\it elementary\/} root unipotents. Since the affine ring of
${\Bbb G}_a$ is $\Int[t]$, one just has to take an independent
variable $t$ as a parameter.
\par
However, it is much harder to figure out, what the generic
element of the {\it elementary\/} subgroup $E(\Phi,\underline{\ \ }\ )$
could be. Also, in the course of localisation proof one has
to construct generic elements of the {\it congruence\/}  subgroups
$G(\Phi,\underline{\ \ },\underline{\ \ }\ )$.
\par
In [47] the second named author finds a way to circumvent these
difficulties. Namely, there he manages to produce a universal
coefficient ring for {\it principal\/} congruence subgroups.
In view of the usual direct limit arguments, it suffices to
carry through a version of localisation method. To be more
specific, let us reproduce the precise statement.
\proclaim
{Theorem~7} There exist a commutative ring $A$, a regular
element $s\in A$ and an element $f\in G(\Phi,A,sA)$
with the following property. For any commutative ring $R$,
any regular element $r\in R$ and any $h\in G(\Phi,R,rR)$
there exists a unique ring homomorphism $\phi:A\map R$
such that $\phi(s)=r$ and $G(\phi)(f)=h$.
\endproclaim
This is the main new tool, but there are many further details,
too technical to be reproduced here. We refer the interested
reader to [47].


\heading
\S~11. Completion
\endheading

Another very important idea, which allows to substantially
enhance the scope of localisation methods, is to combine it
with completion. This idea is due to Anthony Bak [5]. It
was applied to unitary groups in the Thesis of the first named
author [17], [18], and to Chevalley groups, in the works by the
first named and the third named authors [22]. Also, in [17]
and [22] we introduced several important simplifications which
further enhanced the applicability of this method.
\par
Let $s\in R$. Recall that the $s$-completion $\widehat R_{s}$
of the ring $R$ is usually defined as the following inverse
limit:
$$ \widehat R_{s}=\varprojlim R/s^nR,\quad n\in{\Bbb N}. $$
\noindent
However, this definition is not quite compatible with our
purposes. Namely, as always, to control zero divisors, we have
to reduce to Noetherian rings first. Howe\-ver, if $R=\varinjlim
R_i$ is a direct limit of Noetherian rings, the canonical
homomorphism $\varinjlim (\widehat R_i)_{s}\map\widehat R_{s}$ is
in general neither surjective, nor injective.
\par
This forces us to modify the definition of completion as follows:
$$ \widetilde R_{s}=\varinjlim (\widehat R_i)_{s}, $$
\noindent
where the limit is taken over all finitely generated subrings
$R_i$ of $R$ which contain $s$. Let us denote by $\widetilde
F_s$ the canonical map $R\map \widetilde R_{s}$. For the
case, where $R$ is Noetherian $\widetilde F_{s}=\widehat F_{s}$
coincides with the inverse limit of reduction homomorphisms
$\pi_{s^n}:R\map R/s^nR$
\par
Let $R$ be a commutative ring, $\Phi$ be an
irreducible root system of rank $\geq 2$. Define
$$ \aligned
&G(\Phi,R,s^{-1})=\Ker\big(G(\Phi,R)\map
G(\Phi,R_s)/E(\Phi,R_s)\big),\\
\noalign{\vskip 5truept}
&G(\Phi,R,\widehat s)=\Ker\big(G(\Phi,R)\map
G(\Phi,\tilde R_{(s)})/E(\Phi,\tilde R_{(s)})\big).\\
\endaligned $$
\par
The following theorem embodies the gist of
{\bf localisation-completion} method. Morally, it tells that the
commutator of something, that becomes elementary under
localisation, with another something, that becomes elementary
under completion, is indeed elementary. Clearly, this is something
more powerful than just normality of the elementary subgroup.
It goes in the direction of proving that the commutator of
two {\it arbitrary\/} matrices is elementary. Of course, in general
this is not the case, but, as we shall see in the next section,
for finite dimensional rings it is a very near miss.
\proclaim
{Theorem 8C} Let $R$ be a commutative ring, $\Phi$ be an irreducible
root system of rank $\geq 2$. Then
$$ [G(\Phi,R,s^{-1}),G(\Phi,R,\widehat s)]\le E(\Phi,R). $$
\endproclaim
Let $R_i$ be the inductive system of all finitely generated
subrings of $R$, containing $s$. Then
$$ \aligned
&G(\Phi,R,s^{-1})=\varinjlim G(\Phi,R_i,s^{-1}), \\
\noalign{\vskip 5truept}
&G(\Phi,R,\widehat s)=\varinjlim G(\Phi,R_i,\widehat s),\\
\endaligned $$
\noindent
which again reduces the proof to the case, where $R$ is Noetherian.
\par
Let $x\in G(\Phi,R,s^{-1})$ and $y\in G(\Phi,R,\widehat s)$.
By definition, the condition on $x$ means that
$F_s(x)\in E^K(\Phi,a/s^k)$ for some $k$ and $K$. On the other
hand, the condition on $y$ means that
$\pi_{s^n}(y)\in E(\Phi,R/s^nR)$ for {\it all\/} $n$, or, what is
the same, $y\in E(\Phi,R)G(\Phi,R,s^nR)$.
\par
In other words, for {\it any\/} $n$ we can present $y$ as a
product $y=uz$, $u\in E(\Phi,R)$ and $z\in G(\Phi,R,s^nR)$. Thus,
$$ [x,y]=[x,uz]=[x,u]\cdot{}^u[x,z]. $$
\noindent
The first commutator belongs to $E(\Phi,R)$ together with $u$
since $E(\Phi,R)$ is normal. As for the second commutator,
choosing a very large $n$, and applying Theorem 2, we get
$F_s([x,z])\in E(\Phi,F_s(s^qR))$, for a large $q$. On the other
hand, since $G(\Phi,R,s^qR)$ is normal, $[x,z]\in G(\Phi,R,s^qR)$.
Now the usual argument based on the injectivity of the reduction
homomorphism convinces us that $[x,z]\in E(\Phi,s^qR)$.
\par
Similar result holds also in the unitary setting. Let
$(A,\Lambda)$ be a form ring, which is module finite over a
commutative ring $R$. Take $s\in R_0$ and define
$$ \aligned
&U(2n,A,\Lambda,s^{-1})=\Ker\Big(U(2n,A,\Lambda)\map
U(2n,A_s,\Lambda_s)/\EU(2n,A_s,\Lambda_s)\Big),\\
\noalign{\vskip 5truept}
&U(2n,A,\Lambda,\hat s)=\Ker\Big(U(2n,A,\Lambda)\map
U\big(2n,\tilde{(A,\Lambda)}_{(s)}\big)/
E\big(2n,\tilde{(A,\Lambda)}_{(s)}\big)\Big).\\
\endaligned $$
\par
One of the main results of the Thesis by the first named author
[17], [18] can be now stated as follows.
\proclaim
{Theorem~8B} Let $(A,\Lambda)$ be a module finite form ring
over a commutative ring $R$, and let $s\in R_0$. Then
$$ [U(2n,A,\Lambda,s^{-1}),U(2n,A,\Lambda,\hat s)]\le\EU(2n,A,\Lambda). $$
\endproclaim


\heading
\S~12. Nilpotent filtrations
\endheading

As another illustration of the power of our methods, let us state
some important results obtained by the {\bf localisation-completion}
method, as developed in [5], [17], [18], [22], [7]. The following
theorems imply, in particular, nilpotency of {\it relative\/} $\K_1$'s.
\proclaim
{Theorem~9B} Let $(A,\Lambda)$ be a form ring which is module finite
over a commutative ring $R$ of finite Bass---Serre dimension $\delta(R)$,
and let $(I,\Gamma)$ be a form ideal of $(A,\Lambda)$. Then for any
$n\ge 3$ the quotient $U(2n,I,\Gamma)/\EU(2n,I,\Gamma)$ is nilpotent
by abelian of nilpotent class at most $\delta(R)+1$.
\endproclaim
\proclaim
{Theorem~9C} Let $\Phi$ be a reduced irreducible root system of
rank $\ge 2$, let $R$ be a commutative ring of finite
Bass---Serre dimension $\delta(R)$, and let $I\trianglelefteq R$
be its ideal. Then for any Chevalley group $G(\Phi,R)$ of type $\Phi$
over $R$ the quotient $G(\Phi,R,I)/E(\Phi,R,I)$ is nilpotent by abelian
of nilpotent class at most $\delta(R)+1$.
\endproclaim
In fact, in [7] we prove something much more powerful. Namely,
without any finiteness assumptions on ground rings, we construct
{\bf nilpotent filtrations} of congruence subgroups. For rings of
finite Bass---Serre dimension these filtrations are indeed finite.
\par
To state the precise form of these results, we have to recall
definitions of certain {\bf higher elementary subgroups}, which
play the same role for unitary groups, and for Chevalley groups,
as Bak's {\bf very special linear groups} do in the linear case.
\par
Let $(A,\Lambda)$ be a module finite form ring over a commutative ring
$R$ and let $(I,\Gamma)$ be a form ideal in $(A,\Lambda)$. Define
$$ S^d U(2n,I,\Gamma)=\bigcap_{\phi}
\Ker\Big(U(2n,I,\Gamma)\map U(2n,I',\Gamma')/\EU(2n,I',\Gamma')\Big), $$
\noindent
where the intersection is taken over all homomorphisms
$A\map A'$ of rings with involution, $A'$ is module finite over a
commutative ring $R'$ of Bass---Serre dimension $\delta(R')\le d$,
$\Lambda'$ is the form parameter of $A'$ generated by $\phi(\Lambda)$,
$I'$ is the involution invariant ideal of $A'$ generated by $\phi(I)$,
and, finally, $\Gamma'$ is the relative form parameter of
level $I'$ of the form ring $(A',\Lambda')$, generated by
$\phi(\Gamma)$.
\proclaim
{Theorem~10B} Let $(A,\Lambda)$ be a module finite form ring over
a commutative ring $R$.
Further, let $(I,\Gamma)$ be a form ideal of $(A,\Lambda)$
and $n\geq 3$. Then
\smallskip
$\bullet$ Each $S^dU(2n,I,\Gamma)$ is a normal subgroup of
$\GU(2n,R,\Lambda)$.
\smallskip
$\bullet$ The sequence
$$ S^0U(2n,I,\Gamma)\ge S^1U(2n,I,\Gamma)\ge S^2U(2n,I,\Gamma)\ge\cdots $$
\noindent
is a descending $S^0U(2n,R,\Lambda)$-central series.
\smallskip
$\bullet$ The conjugation action of $\GU(2n,R,\Lambda)$ on
$U(2n,I,\Gamma)/S^0U(2n,I,\Gamma)$ is trivial.
\smallskip
$\bullet$ If Bass---Serre dimension of $R$ is finite,
$\delta(R)<\infty$, then
$$ S^dU(2n,R,I)=\EU(2n,R,I), $$
\noindent
whenever $d\ge\delta(R)$.
\endproclaim
Next, we do the same for Chevalley groups. Let $R$ be a commutative
ring and let $I\trianglelefteq R$ an ideal of $R$. Define
$$ S^d G(\Phi,R,I)=\bigcap_{\phi}
\Ker\big(G(\Phi,R,I)\map G(\Phi,A,\phi(I)A)/E(\Phi,A,\phi(I)A)\big), $$
\noindent
where the intersection is taken over all homomorphisms
$\phi:R\map A$ to rings of Bass---Serre dimension $\delta(A)\le d$.
As usual, we set $S^d G(\Phi,R)=S^d G(\Phi,R,R)$.
\proclaim
{Theorem~10C} Let $\Phi$ be an irreducible root system of rank $\ge 2$,
let $R$ be a commu\-ta\-tive ring, and let $I$ be an ideal of $R$. Then
\smallskip
$\bullet$ Each $S^dG(\Phi,R,I)$ is a normal subgroup of $G(\Phi,R)$
\smallskip
$\bullet$ The sequence
$$ S^0G(\Phi,R,I)\ge S^1G(\Phi,R,I)\ge S^2G(\Phi,R,I)\ge\cdots $$
\noindent
is a descending $S^0G(\Phi,A)$-central series.
\smallskip
$\bullet$ The conjugation action of $G(\Phi,R)$ on
$G(\Phi,R,I)/S^0G(\Phi,R,I)$ is trivial.
\smallskip
$\bullet$ If Bass---Serre dimension of $R$ is finite,
$\delta(R)<\infty$, then
$$ S^dG(\Phi,R,I)=E(\Phi,R,I), $$
\noindent
whenever $d\ge\delta(R)$.
\endproclaim


\heading
\S~13. Where next?
\endheading

In conclusion, we list some unsolved problems related to the
results of the present paper. We have preliminary results in
some of these directions, and intend to address them in
subsequent publications.
\proclaim
{Problem~2} Obtain explicit length estimates in the relative
conjugation calculus and commutator calculus.
\endproclaim
\proclaim
{Problem~3} Obtain explicit length estimates in the universal
localisation.
\endproclaim
\proclaim
{Problem~4} Develop versions of universal localisation in the
non-algebraic setting, in particular, for unitary groups.
\endproclaim
Another important challenge is to improve rank bounds in
the commutator calculus for the unitary groups. For the
condition below see [8], [9].
\proclaim
{Problem~5} Develop conjugation calculus and commutator calculus
in the group\/ $\GU(4,R,\Lambda)$, provided\/
$\Lambda R+R\Lambda=\Lambda$.
\endproclaim
\proclaim
{Problem~6}
Prove relative commutator formulae for the group
$\GU(4,R,\Lambda)$, provi\-ded $\Lambda R+R\Lambda=\Lambda$.
\endproclaim
Another important problem is the description of {\it subnormal\/}
subgroups of $G(R)$. For the case of $\GL(n,R)$ this problem has
a fully satisfactory answer, due to the works by John Wilson,
Leonid Vaserstein, and others, see in particular [4], [33], [57],
[61], [63].
\par
For unitary groups, there are works by Gerhard Habdank, the fourth
author, and You Hong, see, in particular, [14], [15], [71] -- [74].
But there are still a number of loose ends.
\proclaim
{Problem~7} Give localisation proofs for the description of
subgroups of the unitary group $\GU(2n,R,\Lambda)$, normalised
by the relative elementary subgroup $\EU(\Phi,I,\Gamma)$,
for a form ideal $(I,\Gamma)$.
\endproclaim
\proclaim
{Problem~8} Using relative localisation, describe subgroups of
a Chevalley group $G(\Phi,R)$, normalised by the relative
elementary subgroup $E(\Phi,R,I)$, for an ideal $I\unlhd R$.
\endproclaim
It would be extremely challenging to fully relativise results
concerning nilpotent filtration.
\proclaim
{Problem~9} Let\/ $R$ be a ring of finite Bass---Serre dimension\/
$\delta(R)=d<\infty$, and let\/ $(I_i,\Gamma_i)$, $1\le i\le m$, be
form ideals of\/ $(R,\Lambda)$. Prove that for any\/ $m>d$ one has
$$ \multline
[[\ldots[G(\Phi,R,I_1),G(\Phi,R,I_2)],\ldots],G(\Phi,R,I_m)]=\\
[[\ldots[E(\Phi,R,I_1),E(\Phi,R,I_2)],\ldots],E(\Phi,R,I_m)].
\endmultline $$
\endproclaim
\proclaim
{Problem~10} Let\/ $R$ be a ring of finite Bass---Serre dimension\/
$\delta(R)=d<\infty$, and let\/ $(I_i,\Gamma_i)$, $1\le i\le m$, be
form ideals of\/ $(R,\Lambda)$. Prove that for any\/ $m>d$ one has
$$ \multline
[[\ldots[\GU(2n,I_1,\Gamma_1),\GU(2n,I_2,\Gamma_2)],\ldots],\GU(2n,I_m,\Gamma_m)]=\\
[[\ldots[\EU(2n,I_1,\Gamma_1),\EU(2n,I_2,\Gamma_2)],\ldots],\EU(2n,I_m,\Gamma_m)].
\endmultline $$
\endproclaim
The following two problems are in fact not individual clear cut
problems, but rather huge research projects.
\proclaim
{Problem~11} Generalise results of the present paper to odd unitary
groups.
\endproclaim
\proclaim
{Problem~12} Obtain results similar to those of the present paper for
{\rm[}groups of points of{\rm]} isotropic reductive groups.
\endproclaim
In the first one of these settings there are foundational works by
Victor Petrov [42] -- [44], while in the second one there are papers
by Victor Petrov, Anastasia Stavrova, and Alexander Luzgarev [45],
[36], with
versions of Quillen---Suslin lem\-ma. But that's about it. Most of
the conjugation calculus and the commutator calculus, including
relative results, explicit estimates, etc., have to be developed
from scratch.


\frenchspacing
\par
\medskip
\widestnumber\no{100}
\newcount\refno
\refno=1

\vfill\eject


\bigskip
Hazrat R., Stepanov A.~V., Vavilov N.~A., Zhang Z.,
The yoga of commutators.
\smallskip
In the present paper we discuss some recent versions of
localisation methods for calculations in the groups of
points of algebraic-like and classical-like groups. Namely,
we describe relative localisation, universal localisation, and
enhanced ver\-sions of localisation-completion. Apart from the
general strategic description of these methods, we state some
typical technical results of the conjugation calculus and the
commutator calculus. Also, we state several recent results obtained
there\-with, such as relative standard commutator formulae, bounded 
width of commuta\-tors, with respect to the elementary generators, 
and nilpotent filtrations of congru\-ence subgroups. Overall, this 
shows that localisation methods can be much more efficient, than 
expected.

\bigskip

Queen's University Belfast, U.K.
\smallskip
{\it E-mail\/}: r.hazrat\@qub.ac.uk, rhazrat\@gmail.com
\medskip

\par
\smallskip
Abdus Salam School of Mathematical Sciences
\par
at the GCU, Lahore, Pakistan
\smallskip
{\it E-mail\/}: stepanov239\@gmail.com
\medskip

\par
Department of Mathematics and Mechanics
\par
St. Petersburg State University\par
St. Petersburg, Russia
\smallskip
{\it E-mail\/}: nikolai-vavilov\@yandex.ru
\medskip

Beijing Institute of Technology, China
\smallskip
{\it E-mail\/}: zuhong\@gmail.com


\bye